# Perturbation Iteration Transform Method for the Solution of Newell-Whitehead-Segel Model Equations

**Grace O. Akinlabi and Sunday O. Edeki**

*Department of Mathematics, Covenant University, Canaanland, Ota, Nigeria*

**Abstract:** In this study, a computational method referred to as Perturbation Iteration Transform Method (PITM), which is a combination of the conventional Laplace Transform Method (LTM) and the Perturbation Iteration Algorithm (PIA) is applied for the solution of Newell-Whitehead-Segel Equations (NWSEs). Three unique examples are considered and the results obtained are compared with their exact solutions graphically. Also, the results agree with those obtained via other semi-analytical methods viz: New Iterative Method and Adomian Decomposition Method. This present method proves to be very efficient and reliable. *Mathematica* 10 is used for all the computations in this study.

## Introduction

The numerical and analytical approximations of Partial Differential Equation (PDE) problems have always been an active field of study in Physics, Mathematics and Engineering. Many researchers have proposed several approaches to solve different PDE problems. For instance, He (2005) solved some wave equations with the Homotopy Perturbation Method (HPM). Akinlabi and Edeki (2016) also solved initial-value wave-like models using the modified Differential Transform Method (DTM). Likewise, Edeki *et al.* (2016) considered the numerical and the analytic solutions of time-fractional linear Schrödinger equations.

In this study, we are concerned with the solution of the Newell-Whitehead-Segel Equations (NWSEs) expressed as follows:

$$\mu_t(x,t) = r\mu_{xx}(x,t) + s\mu(x,t) - t\mu^n(x,t) \qquad (1)$$

Subject to:

$$\mu(x,0) = \varphi(x) \qquad (2)$$

where, $r,s,t \in \mathbb{R}$ with $r>0$ and $n \in \mathbb{Z}^+$.

Equation (1) was derived by Newell and Whitehead (1969; Segel, 1969) and it has been used in modeling various forms of problems that arise from fluid mechanics. It has applications in Chemical, Bio-Engineering and Mechanical Engineering, etc. This equation has been applied to a number of problems. An example of such is in the description of traveling waves by Malomed (1996). Several authors have proposed different methods of solving this equation in the past years. Saravanan and Magesh (2013) solved two nonlinear NWS equations with both the reduced DTM and the Adomian Decomposition Method (ADM). A comparative test was carried out between these two methods where it was shown that the reduced DTM requires less computational work than the ADM. Prakash and Kumar (2016) used the He's Variational Iteration Method to solve non-linear NWSEs, Macías-Díaz and Ruiz-Ramírez (2011) considered the generalised NWSEs using a non-standard symmetry-preserving method. Others are Aasaraai (2011; Ezzati and Shakibi, 2011; Nourazar *et al.*, 2011) and so on.

The idea of using Perturbation Iteration Transform Method (PITM) to solve PDE problems was first pioneered by (Khalid *et al.*, 2016), where the method was used to solve Klein-Gordon equations. The basic idea of this approach is that the PIA and the Laplace Transform (LT) Method are combined to approximate models arising from different fields.

The remaining part of the paper is structured as follows: In section 2 and 3, we review the PIA and the PITM respectively. Section 4 is on the application of the PITM to three cases of NWS equations to show its efficiency. The concluding remark is given in section 5.

## Perturbation Iteration Algorithm (Aksoy and Pakdemirli, 2010)

In this section, we illustrate how the PIA works. If we develop a perturbation algorithm by considering one correction term in the perturbed expansion and the correction terms of the first derivatives in the Taylor series expansion, say $n = 1$ and $m = 1$. The algorithm will be referred to as PIA(1,1).

Now, let us consider a second order differential equation:

$$F(\dot{\mu}, \mu'', \mu, \varepsilon) = 0 \qquad (3)$$

where, $\mu = \mu(x,t)$, $\dot{\mu} = \dfrac{\partial \mu}{\partial t}$, $\mu'' = \dfrac{\partial^2 \mu}{\partial x^2}$ and $\varepsilon$ is the newly introduced perturbed parameter.

And if we consider only one correction term in the expansion:

$$\mu_{n+1} = \mu_n + \varepsilon (\mu_c)_n \qquad (4)$$

Thus, substituting (4) in (3) and expanding such in Taylor series with first derivatives will give:

$$\begin{cases} F(\dot{\mu}, \mu'', \mu, 0) + F_{\dot{\mu}}(\dot{\mu}, \mu'', \mu, 0)\varepsilon(\dot{\mu}_c)_n \\ + F_{\mu''}(\dot{\mu}, \mu'', \mu, 0)\varepsilon(\mu''_c)_n \\ + F_{\mu}(\dot{\mu}, \mu'', \mu, 0)\varepsilon(\mu_c)_n + F_{\varepsilon}(\dot{\mu}, \mu'', \mu, 0)\varepsilon = 0 \end{cases} \qquad (5)$$

where, $\mu = \mu(x,t)$, $F_{\dot{\mu}} = \dfrac{\partial F}{\partial \dot{\mu}}$, $F_{\mu''} = \dfrac{\partial F}{\partial \mu''}$, $F_{\mu} = \dfrac{\partial F}{\partial \mu}$, $F_{\varepsilon} = \dfrac{\partial F}{\partial \varepsilon}$ and $\varepsilon$ the perturbation parameter to be evaluated at $\varepsilon = 0$.

Reorganizing (5), we have:

$$(\dot{\mu}_c)_n + \dfrac{F_{\mu''}}{F_{\dot{\mu}}}(\mu''_c)_n = -\dfrac{F_{\varepsilon} + \dfrac{F}{\varepsilon}}{F_{\dot{\mu}}} - \dfrac{F_{\mu}}{F_{\dot{\mu}}}(\mu_c)_n \qquad (6)$$

With a guessed value, $u_0$, the term $(\mu_c)_0$ is obtained from (6) and then put in (4) to evaluate $\mu_1$. We continue the iterative process using (6) and (4) till the result(s) are satisfied.

## Perturbation Iteration Transform Method (Khalid *et al*., 2016)

To demonstrate the main idea of this method, we consider a PDE with boundary conditions of the form:

$$A\mu(x,t) + B\mu(x,t) + C\mu(x,t) + D\mu(x,t) = h(x,t) \qquad (7)$$

with the initial condition:

$$\mu(x,0) = \varphi(x) \qquad (8)$$

where, $A = \dfrac{\partial}{\partial t}$ is the first order linear differential operator, $B = \dfrac{\partial^2}{\partial x^2}$ is the second order linear differential operator, $C,D$ are the linear and nonlinear terms and $h(x,t)$ is the source term.

Taking the LT of both sides of (7) gives:

$$L[A\mu(x,t)] + L[B\mu(x,t)] + L[C\mu(x,t)] \\ + L[D\mu(x,t)] = L[h(x,t)] \qquad (9)$$

Which on using the differential property of LT yield:

$$L[A\mu(x,t)] = \dfrac{\varphi(x)}{s} + \dfrac{1}{s}L[h(x,t)] \\ - \dfrac{1}{s}L[B\mu(x,t) + C\mu(x,t) + D\mu(x,t)] \qquad (10)$$

Applying the Inverse LT to both sides of (10) gives:

$$\mu(x,t) = H(x,t) \\ - L^{-1}\left[\dfrac{1}{s}L[B\mu(x,t) + C\mu(x,t) + D\mu(x,t)]\right] \qquad (11)$$

where, $H(x,t)$ is the term gotten from the imposed initial condition with the source term.

So, by using the PITM (11) becomes:

$$\mu(x,t) - H(x,t) + \mu_c(x,t)\varepsilon \\ + L^{-1}\left[\dfrac{1}{s}L[B\mu(x,t) + C\mu(x,t) + D\mu(x,t)]\right]\varepsilon = 0 \qquad (12)$$

Thus:

$$\mu_c(x,t) = \dfrac{H(x,t) - \mu(x,t)}{\varepsilon} \\ - L^{-1}\left[\dfrac{1}{s}L[B\mu(x,t) + C\mu(x,t) + D\mu(x,t)]\right] \qquad (13)$$

This is the combined form of the LTM and the PIA. The initial point, $(\mu_c)_0$ is then obtained from (13) and then substituted in (4) to obtain $\mu_1$. The process is iteratively repeated for $\mu_2$, $\mu_3$ and so on. The approximate solution is thus obtained by:

$$\mu(x,t) = \mu_0(x,t) + \mu_1(x,t) + \mu_2(x,t) + \mu_3(x,t) + \cdots \qquad (14)$$

## Illustrative and Numerical Examples

Here, the proposed method is applied to the following problems.

*Cases I, II and III*

*Case I:* Consider the Newell-Whitehead-Segel equation (Patade and Bhalekar, 2015):

$$\mu_t(x,t) = 5\mu_{xx}(x,t) + 2\mu(x,t) + \mu^2(x,t) \quad (15)$$

Subject to:

$$\mu(x,0) = \lambda \quad (16)$$

With the exact solution:

$$\frac{2e^{2t}\lambda}{2+(1-e^{2t})\lambda} \quad (17)$$

*Solution Procedure to Case I*

Taking the LT of both sides of (15), we get:

$$L[\mu(x,t)] = \frac{\lambda}{s} + \frac{1}{s}L[5\mu_{xx}(x,t) + 2\mu(x,t) + \mu^2(x,t)] \quad (18)$$

Applying the ILT to both sides of (18) gives:

$$\mu(x,t) = \lambda + L^{-1}\left[\frac{1}{s}L[5\mu_{xx}(x,t) + 2\mu(x,t) + \mu^2(x,t)]\right] \quad (19)$$

Now, by the PITM (19) becomes:

$$\mu(x,t) - \lambda + \mu_c(x,t)\varepsilon$$
$$-L^{-1}\left[\frac{1}{s}L[5\mu_{xx}(x,t) + 2\mu(x,t) + \mu^2(x,t)]\right]\varepsilon = 0 \quad (20)$$

Thus:

$$\mu_c(x,t) = \frac{-\mu(x,t)+\lambda}{\varepsilon}$$
$$+L^{-1}\left[\frac{1}{s}L[5\mu_{xx}(x,t) + 2\mu(x,t) + \mu^2(x,t)]\right] \quad (21)$$

This implies that:

$$\mu_0(x,t) = \lambda$$
$$\mu_1(x,t) = 2t\lambda + t\lambda^2$$
$$\mu_2(x,t) = 2t^2\lambda + t^2\lambda^2 + \frac{4t^3\lambda^2}{3} + \frac{4t^3\lambda^3}{3} + \frac{t^3\lambda^4}{3}$$
$$\mu_3(x,t) = \frac{4t^3\lambda}{3} + \frac{2t^3\lambda^2}{3} + \frac{2t^4\lambda^2}{3} + \frac{4t^5\lambda^2}{5} + \frac{2t^4\lambda^2}{3} + \frac{4t^5\lambda^3}{5}$$
$$+\frac{8t^6\lambda^3}{8} + \frac{t^4\lambda^4}{6} + \frac{t^5\lambda^4}{5} + \frac{4t^6\lambda^4}{3} + \frac{16t^7\lambda^4}{63} + \frac{2t^6\lambda^5}{3}$$
$$+\frac{32t^7\lambda^5}{63} + \frac{t^6\lambda^6}{9} + \frac{8t^7\lambda^6}{21} + \frac{8t^7\lambda^7}{63} + \frac{t^7\lambda^8}{63}$$
$$\vdots$$

Therefore, the solution $\mu(x,t)$ is given by:

$$\mu(x,t) = \mu_0(x,t) + \mu_1(x,t) + \mu_2(x,t) + \mu_3(x,t) + \cdots$$
$$= \lambda + 2t\lambda + t\lambda^2 + 2t^2\lambda + t^2\lambda^2 + \frac{4t^3\lambda^2}{3} + \frac{4t^3\lambda^3}{3}$$
$$+\frac{t^3\lambda^4}{3} + \frac{4t^3\lambda}{3} + \frac{2t^3\lambda^2}{3} + \frac{2t^4\lambda^2}{3} + \frac{4t^5\lambda^2}{5} + \frac{2t^4\lambda^2}{3} \quad (22)$$
$$+\frac{4t^5\lambda^3}{5} + \frac{8t^6\lambda^3}{8} + \frac{t^4\lambda^4}{6} + \frac{t^5\lambda^4}{5} + \frac{4t^6\lambda^4}{3} + \frac{16t^7\lambda^4}{63}$$
$$+\frac{2t^6\lambda^5}{3} + \frac{32t^7\lambda^5}{63} + \frac{t^6\lambda^6}{9} + \frac{8t^7\lambda^6}{21} + \frac{8t^7\lambda^7}{63} + \frac{t^7\lambda^8}{63} + \cdots$$

*Case II:* Consider the Newell-Whitehead-Segel equation (Patade and Bhalekar, 2015):

$$\mu_t(x,t) = \mu_{xx}(x,t) + 2\mu(x,t) - 3\mu^2(x,t) \quad (23)$$

Subject to:

$$\mu(x,0) = \lambda \quad (24)$$

With the exact solution:

$$\frac{2e^{2t}\lambda}{-2+3(1-e^{2t})\lambda} \quad (25)$$

*Solution Procedure to Case II*

Taking the LT of both sides of (23), we get:

$$L[\mu(x,t)] = \frac{\lambda}{s} + \frac{1}{s}L[\mu_{xx}(x,t) + 2\mu(x,t) - 3\mu^2(x,t)] \quad (26)$$

Applying the ILT to both sides of (26) gives:

$$\mu(x,t) = \lambda + L^{-1}\left[\frac{1}{s}L[\mu_{xx}(x,t) + 2\mu(x,t) - 3\mu^2(x,t)]\right] \quad (27)$$

Now, by the PITM (27) becomes:

$$\mu(x,t) - \lambda + \mu_c(x,t)\varepsilon$$
$$-L^{-1}\left[\frac{1}{s}L[\mu_{xx}(x,t) + 2\mu(x,t) - 3\mu^2(x,t)]\right]\varepsilon = 0 \quad (28)$$

Thus:

$$\mu_c(x,t) = \frac{-\mu(x,t)+\lambda}{\varepsilon}$$
$$+L^{-1}\left[\frac{1}{s}L[\mu_{xx}(x,t) + 2\mu(x,t) - 3\mu^2(x,t)]\right] \quad (29)$$

This implies that:

$\mu_0(x,t) = \lambda$

$\mu_1(x,t) = 2t\lambda - 3t\lambda^2$

$\mu_2(x,t) = 2t^2\lambda - 3t^2\lambda^2 - 4t^3\lambda^2 + 12t^3\lambda^3 - 9t^3\lambda^4$

$\mu_3(x,t) = \dfrac{4t^3\lambda}{3} - 2t^3\lambda^2 - 2t^4\lambda^2 - \dfrac{12t^5\lambda^2}{5} + 6t^4\lambda^3 + \dfrac{36t^5\lambda^3}{5}$

$+ 8t^6\lambda^3 - \dfrac{9t^4\lambda^4}{2} - \dfrac{27t^5\lambda^4}{5} - 36t^6\lambda^4 - \dfrac{48t^7\lambda^4}{7} + 54t^6\lambda^5$

$+ \dfrac{288t^7\lambda^5}{7} - 27t^6\lambda^6 - \dfrac{648t^7\lambda^6}{7} + \dfrac{648t^7\lambda^7}{7} - \dfrac{243t^7\lambda^8}{7}$

$\vdots$

Therefore, the solution $\mu(x,t)$ is given by:

$\mu(x,t) = \mu_0(x,t) + \mu_1(x,t) + \mu_2(x,t) + \mu_3(x,t) + \cdots$

$= \lambda + 2t\lambda - 3t\lambda^2 + 2t^2\lambda - 3t^2\lambda^2 - 4t^3\lambda^2 + 12t^3\lambda^3 - 9t^3\lambda^4$

$+ \dfrac{4t^3\lambda}{3} - 2t^3\lambda^2 - 2t^4\lambda^2 - \dfrac{12t^5\lambda^2}{5} + 6t^4\lambda^3 + \dfrac{36t^5\lambda^3}{5}$ (30)

$+ 8t^6\lambda^3 - \dfrac{9t^4\lambda^4}{2} - \dfrac{27t^5\lambda^4}{5} - 36t^6\lambda^4 - \dfrac{48t^7\lambda^4}{7} + 54t^6\lambda^5$

$+ \dfrac{288t^7\lambda^5}{7} - 27t^6\lambda^6 - \dfrac{648t^7\lambda^6}{7} + \dfrac{648t^7\lambda^7}{7} - \dfrac{243t^7\lambda^8}{7} + \cdots$

*Case III:* Consider the Newell-Whitehead equation (Patade and Bhalekar, 2015):

$\mu_t(x,t) = \mu_{xx}(x,t) + 2\mu(x,t) - 3\mu^3(x,t)$ (31)

Subject to:

$\mu(x,0) = \sqrt{\dfrac{2}{3}} \dfrac{e^{2x}}{e^x + e^{2x}}$ (32)

*Solution Procedure to Case III*

Taking the LT of both sides of (31), we get:

$L[\mu(x,t)] = \dfrac{1}{s}\sqrt{\dfrac{2}{3}} \dfrac{e^{2x}}{e^x + e^{2x}}$

$+ \dfrac{1}{s} L[\mu_{xx}(x,t) + 2\mu(x,t) - 3\mu^3(x,t)]$ (33)

Applying the ILT to both sides of (33) gives:

$\mu(x,t) = \sqrt{\dfrac{2}{3}} \dfrac{e^{2x}}{e^x + e^{2x}}$

$+ L^{-1}\left[\dfrac{1}{s} L[\mu_{xx}(x,t) + 2\mu(x,t) - 3\mu^3(x,t)]\right]$ (34)

Now, by the PITM (34) becomes:

$\mu(x,t) - \sqrt{\dfrac{2}{3}} \dfrac{e^{2x}}{e^x + e^{2x}} + \mu_c(x,t)\varepsilon$

$- L^{-1}\left[\dfrac{1}{s} L[\mu_{xx}(x,t) + 2\mu(x,t) - 3\mu^3(x,t)]\right]\varepsilon = 0$ (35)

Thus:

$\mu_c(x,t) = \dfrac{-\mu(x,t) + \sqrt{\dfrac{2}{3}} \dfrac{e^{2x}}{e^x + e^{2x}}}{\varepsilon}$

$+ L^{-1}\left[\dfrac{1}{s} L[\mu_{xx}(x,t) + 2\mu(x,t) - 3\mu^3(x,t)]\right]$ (36)

This implies that:

$\mu_0(x,t) = \sqrt{\dfrac{2}{3}} \dfrac{e^{2x}}{e^x + e^{2x}}$, $\mu_1(x,t) = -\dfrac{2\sqrt{\tfrac{2}{3}} e^{6x}}{(e^x + e^{2x})^3} t$

$+ \dfrac{2\sqrt{6} e^{2x}}{e^x + e^{2x}} t - \dfrac{4\sqrt{\tfrac{2}{3}} e^{2x}(e^x + 2e^{2x})}{(e^x + e^{2x})^2} t$

$+ \sqrt{\dfrac{2}{3}} e^{2x}\left(\dfrac{2(e^x + 2e^{2x})^2}{(e^x + e^{2x})^3} - \dfrac{e^x + 4e^{2x}}{(e^x + e^{2x})^2}\right) t$

$\mu_2(x,t) = \dfrac{3\sqrt{\tfrac{3}{2}} e^x (1 + e^{2x}) t^2}{(1 + e^x)^4} - \dfrac{9\sqrt{\tfrac{3}{2}} e^{3x} t^4}{(1 + e^x)^6}$

$\mu_3(x,t) = \dfrac{\sqrt{\tfrac{3}{2}} e^x (3 - 6e^x + 22e^{2x} - 6e^{3x} + 3e^{4x}) t^3}{(1 + e^x)^6}$

$- \dfrac{9\sqrt{\tfrac{3}{2}} e^{3x}(11 - 20e^x + 11e^{2x}) t^5}{5(1 + e^x)^8} - \dfrac{243\sqrt{\tfrac{3}{2}} e^{3x}(1 + e^{2x})^3 t^7}{14(1 + e^x)^{12}}$

$+ \dfrac{243\sqrt{\tfrac{3}{2}} e^{5x}(1 + e^{2x})^2 t^9}{2(1 + e^x)^{14}} - \dfrac{6561\sqrt{\tfrac{3}{2}} e^{7x}(1 + e^{2x}) t^{11}}{22(1 + e^x)^{16}}$

$+ \dfrac{6561\sqrt{\tfrac{3}{2}} e^{9x} t^{13}}{26(1 + e^x)^{18}}$

$\vdots$

Therefore, the solution $\mu(x,t)$ is given by:

$\mu(x,t) = \mu_0(x,t) + \mu_1(x,t) + \mu_2(x,t) + \mu_3(x,t) + \cdots$

$= \sqrt{\dfrac{2}{3}} \dfrac{e^{2x}}{e^x + e^{2x}} - \dfrac{2\sqrt{\tfrac{2}{3}} e^{6x}}{(e^x + e^{2x})^3} t + \dfrac{2\sqrt{6} e^{2x}}{e^x + e^{2x}} t - \dfrac{4\sqrt{\tfrac{2}{3}} e^{2x}(e^x + 2e^{2x})}{(e^x + e^{2x})^2} t$

$+ \sqrt{\dfrac{2}{3}} e^{2x}\left(\dfrac{2(e^x + 2e^{2x})^2}{(e^x + e^{2x})^3} - \dfrac{e^x + 4e^{2x}}{(e^x + e^{2x})^2}\right) t + \dfrac{3\sqrt{\tfrac{3}{2}} e^x (1 + e^{2x}) t^2}{(1 + e^x)^4}$

$- \dfrac{9\sqrt{\tfrac{3}{2}} e^{3x} t^4}{(1 + e^x)^6} + \dfrac{\sqrt{\tfrac{3}{2}} e^x (3 - 6e^x + 22e^{2x} - 6e^{3x} + 3e^{4x}) t^3}{(1 + e^x)^6}$ (37)

$- \dfrac{9\sqrt{\tfrac{3}{2}} e^{3x}(11 - 20e^x + 11e^{2x}) t^5}{5(1 + e^x)^8} - \dfrac{243\sqrt{\tfrac{3}{2}} e^{3x}(1 + e^{2x})^3 t^7}{14(1 + e^x)^{12}}$

$+ \dfrac{243\sqrt{\tfrac{3}{2}} e^{5x}(1 + e^{2x})^2 t^9}{2(1 + e^x)^{14}} - \dfrac{6561\sqrt{\tfrac{3}{2}} e^{7x}(1 + e^{2x}) t^{11}}{22(1 + e^x)^{16}}$

$+ \dfrac{6561\sqrt{\tfrac{3}{2}} e^{9x} t^{13}}{26(1 + e^x)^{18}} + \cdots$

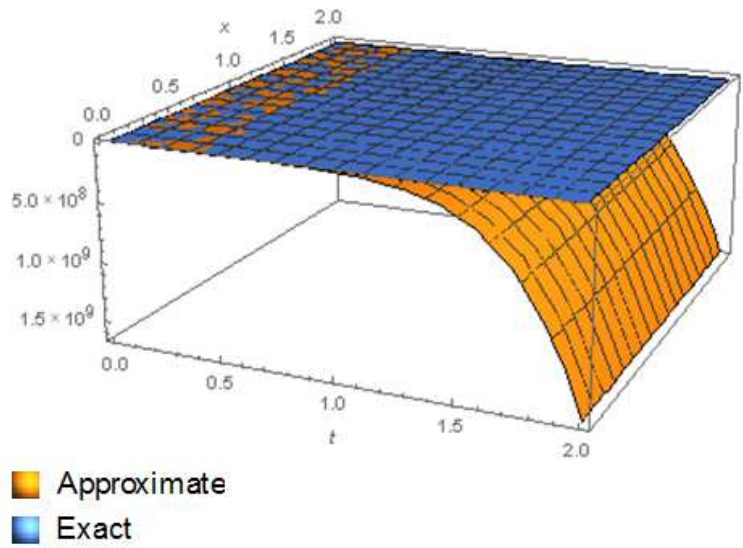

Fig. 1. Solution graph for case I

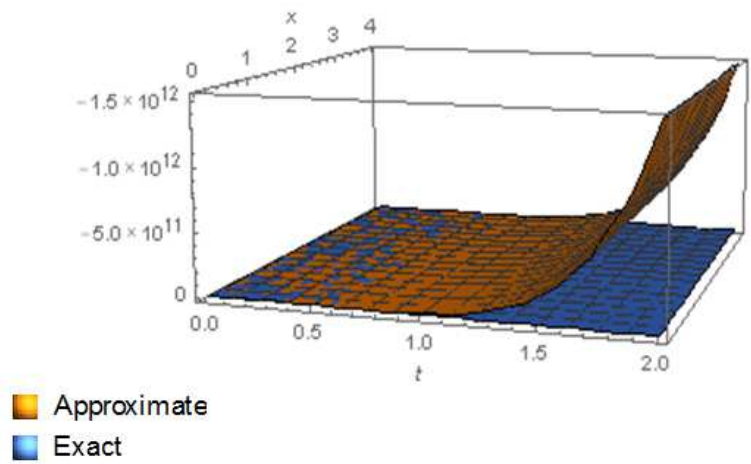

Fig. 2. Solution graph for case II

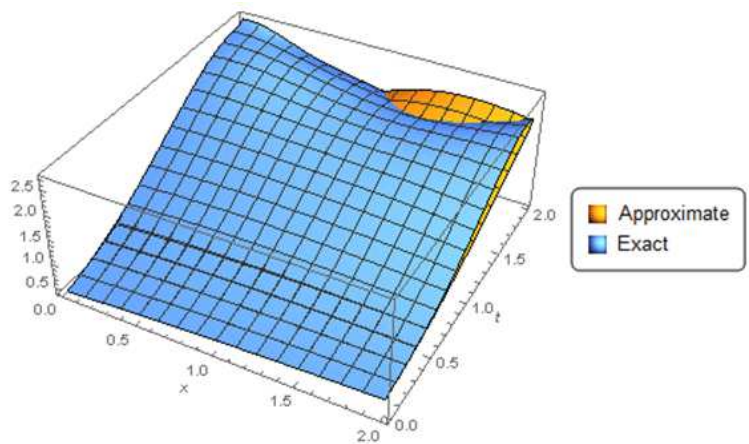

Fig. 3. Solution graph for case III

■■■

Figure 1 is for the exact solution and the PITM (approximate solution) of case I.

Figure 2 is for exact solution and the PITM (approximate solution) of case II.

Figure 3 is for exact solution and the PITM (approximate solution) of case III.

*Discussion of Results*

In this subsection, we present the graphs for the exact solution and the PITM (approximate solution) for the three cases.

## Conclusion

In this study, the solutions (roots) of the Newell-Whitehead-Segel models are gotten using Perturbation Iteration Transform Method as a proposed computational method. The results are obtained with less computational time and are compared graphically with their exact solutions. In addition, these results are in good agreement with those by Patade and Bhalekar (2015) using other semi-analytical method: Adomian Decomposition Method and New Iterative Method. We therefore, propose this method for the solutions of both linear and non-linear PDEs in other aspects of applied sciences.

## Acknowledgement

The authors wish to sincerely thank Covenant University (CUCRID) for the provision of enabling environment and finance. In addition, the authors thank the anonymous referees for their valuable contributions.

## Authors' Contributions

The concerned authors: GOA and SOE contributed significantly to this work. Both authors approved the final manuscript for publication.

## Ethics

The authors declare that there exists no conflict of interest regarding this paper.

## References

Aasaraai, A., 2011. Analytic solution for Newell-whitehead-segel equation by differential transform method. Middle East J. Scientific Res., 10: 270-273.

Akinlabi, G.O. and S.O. Edeki, 2016. On approximate and closed-form solution method for initial-value wave-like models. Int. J. Pure Applied Math., 107: 449-456. DOI: 10.12732/ijpam.v107i2.14

Aksoy, Y. and M. Pakdemirli, 2010. New perturbation-iteration solutions for Bratu-type equations. Comput. Math. Applic., 59: 2802-2808. DOI: 10.1016/j.camwa.2010.01.050

Edeki, S.O., G.O. Akinlabi and S.A. Adeosun, 2016. Analytic and numerical solutions of time-fractional linear schrödinger equation. Commun. Math. Applic., 7: 1-10.

Ezzati, R. and K. Shakibi, 2011. Using adomian's decomposition and multiquadric quasi-interpolation methods for solving Newell-Whitehead equation. Proc. Comput. Sci., 3: 1043-1048. DOI: 10.1016/j.procs.2010.12.171

He, J.H., 2005. Application of homotopy perturbation method to nonlinear wave equations. Chaos Solitons Fractals, 26: 695-700. DOI: 10.1016/j.chaos.2005.03.006

Khalid, M., M. Sultana, F. Zaidi and A. Uroosa, 2016. Solving linear and nonlinear Klein-Gordon equations by new perturbation iteration transform method. TWMS J. Applied Eng. Math., 6: 115-125.

Macías-Díaz, J. and J. Ruiz-Ramírez, 2011. A non-standard symmetry-preserving method to compute bounded solutions of a generalized Newell-Whitehead-Segel equation. Applied Numerical Math., 61: 630-640. DOI: 10.1016/j.apnum.2010.12.008

Malomed, B.A., 1996. The Newell-Whitehead-Segel equation for travelling waves.

Newell, A. and J. Whitehead, 1969. Finite bandwidth, finite amplitude convection. J. Fluid Mechan., 38: 279-303. DOI: 10.1017/S0022112069000176

Nourazar, S.S., M. Soori and A. Nazari-Golshan, 2011. On the exact solution of Newell-Whitehead-Segel equation using the homotopy perturbation method. Australian J. Basic Applied Sci. Fractals, 5: 1400-1411.

Patade, J. and S. Bhalekar, 2015. Approximate analytical solutions of Newell-Whitehead-Segel equation using a new iterative method. World J. Modell. Simulat., 11: 94-103.

Prakash, A. and M. Kumar, 2016. He's variational iteration method for the solution of nonlinear Newell-Whitehead-Segel equation. J. Applied Anal. Comput., 6: 738-748. DOI: 10.11948/2016048

Saravanan, A. and N. Magesh, 2013. A comparison between the reduced differential transform method and the Adomian decomposition method for the Newell-Whitehead-Segel equation. J. Egypt. Math. Society, 21: 259-265. DOI: 10.1016/j.joems.2013.03.004

Segel, L.A., 1969. Distant side-walls cause slow amplitude modulation of cellular convection. J. Fluid Mechan., 38: 203-224. DOI: 10.1017/S0022112069000127